\documentclass[twocolumn]{autart}    

\usepackage{graphicx}          
\usepackage{epstopdf}
\usepackage{fancyhdr}
\usepackage{indentfirst}
\usepackage{amsmath}
\usepackage{amssymb}
\usepackage[longnamesfirst,round]{natbib}
\usepackage{natbib}
\usepackage{color}
\usepackage{latexsym}
\usepackage{amsfonts}
\usepackage{amssymb,amsmath,amsfonts}
\usepackage{amsmath,mathrsfs,bm,url,times}
\usepackage{caption}
\usepackage{subcaption}
\usepackage{textcomp}
\usepackage{epsfig}
\newtheorem{theorem}{Theorem}
\newtheorem{lemma}{Lemma}
\newtheorem{definition}{Definition}
\newtheorem{remark}{Remark}
\DeclareMathOperator{\sat}{sat}
\DeclareMathOperator{\scc}{SCC}
\DeclareMathOperator{\argmax}{argmax}
\DeclareMathOperator{\diag}{diag}

\newenvironment{proof}{\vspace{1ex}\noindent{\bf Proof.}\hspace{0.5em}}
    {\hfill\qed\vspace{1ex}}

\allowdisplaybreaks
\begin{document}

\begin{frontmatter}

\title{Distributed Event-Triggered Control for Global Consensus of Multi-Agent Systems with Input Saturation\thanksref{footnoteinfo}} 

\thanks[footnoteinfo]{This work was supported by the
Knut and Alice Wallenberg Foundation, the  Swedish Foundation for Strategic Research, and the Swedish Research Council. This paper was not presented at any IFAC
meeting. Corresponding author: T.~Yang, Tel: +1 940-891-6876,
Fax: +1 940-891-6881.}

\author[Paestum]{Xinlei Yi}\ead{xinleiy@kth.se},    
\author[Rome]{Tao Yang}\ead{Tao.Yang@unt.edu},               
\author[Paestum]{Junfeng Wu}\ead{junfengw@kth.se}, 
\author[Paestum]{Karl H. Johansson}\ead{kallej@kth.se}

\address[Paestum]{ACCESS Linnaeus Centre, Electrical Engineering, KTH Royal Institute of Technology, 100 44, Stockholm, Sweden}  
\address[Rome]{Department of Electrical Engineering, University of North Texas, Denton, TX 76203 USA}             


\begin{abstract}                          
We consider the global consensus problem for multi-agent systems with input saturation over digraphs. Under a mild connectivity condition that the underlying digraph has a directed spanning tree, we use Lyapunov methods to show that the widely used distributed consensus protocol, which solves the consensus problem for the case without input saturation constraints,
also solves the global consensus problem for the case with input saturation constraints. In order to reduce the overall need of communication and system updates,
we then propose a distributed event-triggered control law. Global consensus is still realized and Zeno behavior is excluded. Numerical simulations are provided to illustrate the effectiveness of the theoretical results.
\end{abstract}

\end{frontmatter}

\section{INTRODUCTION}
\label{sec:intro}
In the past decades, consensus in multi-agent systems has been
widely investigated because of its wide applicability. In the widely used distributed consensus protocol setup,
each agent updates its state based on its own and the states of its
neighbors in such a way that the final states of all agents converge to a
common value. Consensus problem has been studied extensively \citep[e.g.,][]{olfati2004consensus,ren2007information,Liu2011consensus,you2011network} and the references therein.
It is known consensus is achieved if the underlying graph is directed and has a spanning tree.

However, real systems are subject to physical constraints, such as input,
output, digital communication channels, and sensors constraints. These constraints lead to nonlinearity in the closed-loop dynamics. Thus the behavior of each agent is affected and special attention to these constraints needs to be taken in order to understand their influence on the convergence properties. Here we list some representative examples of such constraints. For example, \citet{yang2014global} studies global consensus for discrete-time multi-agent systems with input saturation constraint; \citet{meng2013global} considers the leader-following consensus problem for multi-agent systems subject to input saturation; \citet{lim2016consensus} and \citet{wang2016conditions} investigate necessary and sufficient initial conditions for achieving consensus in the presence of
output saturation; \citet{li2011consensus} shows that the widely used distributed consensus protocol also asymptotically leads to consensus for multi-agent systems with input saturation and directed topologies. 

The widely used distributed consensus protocols (see e.g., \cite{olfati2004consensus}) require continuous information exchange among the agents. For example, each agent needs to continuously broadcast its state to its neighbors.
Therefore, it may be impractical to require continuous communication in physical
applications. On the other hand, motivated by the future trend that agents can be equipped with embedded microprocessors with limited capabilities to transmit and collect data and do simple data processing, event-triggered control is introduced partially to tackle this problem \citep{aastrom1999comparison,tabuada2007event,wang2011event,heemels2012introduction}. The
control in event-triggered control is often piecewise constant between
triggering times. The triggering times are determined implicitly by the event conditions. Event-triggered control for multi-agent systems has been extensively studied by many researchers recently \citep[e.g.,][]{dimarogonas2012distributed,seyboth2013event,meng2013event,fan2013distributed,meng2015periodic,yang2016decentralized,yi2016distributed,yi2016formation,yi2017pull}.
A key challenge in event-triggered control for multi-agent systems is how to design control law, the event threshold to determine the triggering times, and to exclude Zeno behavior. For continuous-time multi-agent systems, Zeno behavior means there are infinite number of triggers in a
finite time interval \cite{johansson1999regularization}. 
Another important problem is how to realize the event-triggered controller in a distributed way.

It is well known that, in almost all real applications, actuators have bounds on their inputs and thus actuator saturation is important to study. However, all the papers mentioned above do not take input saturation into consideration. The consensus problem for multi-agent systems with input saturation and event-triggered controllers is a challenging problem since these constraints lead to nonlinearity in the closed-loop dynamics. Thus the behavior of each agent is affected and special attention to these constraints needs to be taken in order to understand their influence on the convergence properties. To the best of our knowledge, there are only few papers addressed this. \citet{wu2016distributed} proposes the distributed event-triggered
control strategy to achieve consensus for multi-agent systems which is subject to input saturation in an output feedback mechanism. Different from this paper, the underlying graph they used was undirected and they did not consider excluding Zeno behavior. Actually, even for a single agent system with input saturation and the event-triggered controller, the stability problem is also challenging. \cite{kiener2014actuator} addresses the influence of actuator saturation on event-triggered control. \citet{xie2017event} studies the problem of global stabilization
of multiple integrator systems using event-triggered bounded controls.

In this paper, we address the global consensus problem for multi-agent systems with input saturation over digraphs. Specifically,
we first consider the case that the underlying graph is directed and strongly connected. In this case, we show that, under the widely used distributed consensus protocol, multi-agent systems with input saturation achieve consensus. Then, we study the case that the underlying graph is directed and has a directed spanning tree, and show that consensus can still be reached. Finally, in order to reduce  actuation updates and inter-agent communications, we propose an event-triggered law which also leads  to that consensus is achieved.  Our main contributions are twofold: (1) a Lyapunov function which is different from the one in \citet{li2011consensus} is used to prove our results, which facilitates the design of event-triggered control law; (2) event-triggered control law in this paper is distributed in the sense that it does not require any a priori knowledge of global network parameters and it is free from Zeno behavior.

The remainder of this paper is organized as follows. Section \ref{secpreliminaries} introduces the preliminaries and the problem formulation. The main results are stated in Section \ref{secmain} and Section \ref{secmaini}. Simulations are given in Section \ref{secsimulation}. Finally, the paper is concluded in Section \ref{secconclusion}.

\noindent {\bf Notations}: $\|\cdot\|$ represents the Euclidean norm for
vectors or the induced 2-norm for matrices. ${\bf 1}_n$ denotes the column
vector with each component 1 and dimension $n$. $I_n$ is the $n$ dimension identity matrix. $\rho(\cdot)$ stands for
the spectral radius for matrices and $\rho_2(\cdot)$ indicates the minimum
positive eigenvalue for matrices having positive eigenvalues. Given two
symmetric matrices $M,N$, $M>N$ (or $M\ge N$) means $M-N$ is a positive
definite (or positive semi-definite) matrix. The notation $A\otimes B$ denotes the Kronecker product
of matrices (vectors) $A$ and $B$. Given a vector $s=[s_1,\dots,s_n]\in\mathbb{R}^n$, define the component operators $c_l(s)=s_l,l=1,\dots,n$. 

\section{PRELIMINARIES}\label{secpreliminaries}
In this section, we present some definitions from algebraic graph theory \citep{mesbahi2010graph} and the problem formulation.

\subsection{Algebraic Graph Theory}

Let $\mathcal G=(\mathcal V,\mathcal E, \mathcal A)$ denote a  (weighted) directed graph (or digraph) with the set of agents (vertices or nodes) $\mathcal V =\{v_1,\cdots,v_n\}$, the set of links (edges) $\mathcal E\subseteq \mathcal V \times \mathcal V$, and the weighted adjacency matrix
$\mathcal A =(a_{ij})$ with nonnegative adjacency elements $a_{ij}$. A link
of $\mathcal G$ is denoted by $(v_i,v_j)\in \mathcal E$ if there is
a directed link from agent $v_j$ to agent $v_i$ with weight $a_{ij}>0$, i.e. agent $v_j$ can send
information to agent $v_i$ while the opposite direction transmission might not exist or with different weight $a_{ji}$. The adjacency elements associated with the links of the graph
are positive, i.e., $(v_i,v_j)\in \mathcal E\iff a_{ij}>0$. It is assumed that $a_{ii}=0$ for all $i\in\mathcal I$, where $\mathcal I=\{1\dots,n\}$.
The in-degree of agent $v_i$ is defined as $
\deg^{\text{in}}_i=\sum\limits_{j=1}^{n}a_{ij}$.
The degree matrix of digraph $\mathcal G$ is defined as
$D=\diag([\deg^{\text{in}}_1, \cdots, \deg^{\text{in}}_n])$. The Laplacian
matrix associated with the digraph $\mathcal G$ is defined as $L=D-\mathcal
A$. A directed path from agent $v_0$ to agent $v_k$ is a directed graph
with distinct agents $v_0,...,v_k$ and links $e_0,...,e_{k-1}$ such that
$e_i$ is a link directed from $v_i$ to $v_{i+1}$, for all $i<k$.

\begin{definition}
A directed graph $\mathcal G$ is strongly connected if for any two
distinct agents $v_i,v_j$, there exits a directed path from agent $v_i$ to agent
$v_j$.
\end{definition}

By \citet[Theorem 6.2.24]{horn2012matrix}, we know that strongly connectivity of $\mathcal G$  is
equivalent to the irreducibility of the corresponding Laplacian matrix $L$.
\begin{definition}
A directed graph $\mathcal G$ has a directed spanning tree if there exists
one agent $v_{i_{0}}$ such that for any other agent $v_{j}$, there exits a
directed path from $v_{i_{0}}$ to $v_{j}$.
\end{definition}

By Perron-Frobenius theorem \citep{horn2012matrix} (for more details and proof, see
\citet{lu2006new} and \citet{lu2007new}), we have
\begin{lemma}\label{lemmaL}
Suppose $L$ is the Laplacian matrix associated with a digraph $\mathcal G$ that has a spanning tree,
then $rank(L)=n-1$, and zero is an algebraically simple
eigenvalue of $L$. Moreover, if $L$ is irreducible, there is a positive vector
$\xi^{\top}=[\xi_{1},\cdots,\xi_{n}]$ such that $\xi^{\top} L=0$ and
$\sum_{i=1}^{n}\xi_{i}=1$. 
\end{lemma}
The following result in \citet{yi2017pull} is also useful for our analysis later.
\begin{lemma}\label{lemmaU}
Suppose that $L$ is irreducible and $\xi$ is the vector defined in Lemma \ref{lemmaL}.
Let $\Xi=\diag(\xi)$, $U=\Xi-\xi\xi^\top$, and $R=\frac{1}{2}(\Xi L+L^\top\Xi)$. Then $R=\frac{1}{2}(UL+L^\top U)$ and
\begin{align}\label{ULL}
U\ge\frac{\rho_2(U)}{\rho(L^\top L)}L^\top L\ge0~\text{and}~R\ge \frac{\rho_2(R)}{\rho(U)}U\ge0.
\end{align}
\end{lemma}

\subsection{Multi-Agent Systems with Input Saturation}
We consider a set of $n$ agents that are modelled as a single integrator with input saturation:
\begin{align}
\dot{x}_i(t)=\sat_h(u_i(t))
,~i\in\mathcal I, t\ge0,\label{system}
\end{align}
where $x_i(t)\in\mathbb{R}^p$ is the state and $u_i(t)\in\mathbb{R}^p$ is the control input of agent $v_i$, respectively. And $\sat_h(\cdot)$ is the saturation function defined as
\begin{align}\label{satufun}
\sat_h(s)=[\sat_h(s_1),\dots,\sat_h(s_l)]^\top,
\end{align}
where $s=[s_1,\dots,s_l]^\top\in\mathbb{R}^l$ with $l$ is a positive integer and
\begin{align*}
\sat_h(s_i)=
\begin{cases}
h,~&\text{if}~s_i\ge h\\
s_i,~&\text{if}~|s_i|<h\\
-h,~&\text{if}~s_i\le -h
\end{cases},
\end{align*}
with $h$ is a positive constant, and the interval $[-h,h]$ is referred to as the saturation level.

\begin{remark}
For the ease of presentation, in this paper, we focus on the case where all the agents have the same input saturation level. However, the analysis can be extended to the case where the agents have different saturation levels.
\end{remark}

In the  literature, the widely used distributed consensus protocol is
\begin{align}\label{inputc}
u_i(t)=-\sum_{j=1}^{n}L_{ij}x_j(t).
\end{align}

In this paper, we first show that consensus can be achieved even in the presence of input saturation, i.e., consensus is achieved for the multi-agent system (\ref{system}) with the distributed protocol (\ref{inputc}).

The following properties about the saturation function are useful for our analysis.
\begin{lemma}\label{lemma3}
For any constants $a$ and $b$,
\begin{align*}
&\frac{1}{2}a^2\ge\int_{0}^{a}\sat_h(s)ds\ge\frac{1}{2}(\sat_h(a))^2,\\
&(a-b)^2\ge(\sat_h(a)-\sat_h(b))^2.
\end{align*}
\end{lemma}

\begin{lemma}\label{lemmaLx}
Suppose that $L$ is the Laplacian matrix associated with a digraph $\mathcal G$ that has a spanning tree. For $x_1,\dots,x_n\in\mathbb{R}^p$, define $y_i=\sat_h(-\sum_{j=1}^{n}L_{ij}x_j)\in\mathbb{R}^p$. Then $y_1=\cdots=y_n$ if and only if $x_1=\cdots=x_n$.
\end{lemma}
{\bf Proof}: The sufficiency is obvious.  Let's show the necessity.
 Let $z_i=-\sum_{j=1}^{n}L_{ij}x_j$. From $y_1=\cdots=y_n$, we know that for any $l=1,\dots,p$, $c_l(z_i)>0,\forall i\in\mathcal I$, or $c_l(z_i)<0,\forall i\in\mathcal I$, or $c_l(z_i)=0,\forall i\in\mathcal I$.

From Lemma 2 in \citet[Lemma~2]{li2011consensus}, we know that neither $c_l(z_i)>0,\forall i\in\mathcal I$ nor $c_l(z_i)<0,\forall i\in\mathcal I$ holds. Thus $c_l(z_i)=0,\forall i\in\mathcal I$. Then, from Lemma \ref{lemmaL}, we have $c_l(x_i)=c_l(x_j),\forall i,j\in\mathcal I$. Hence $x_1=\cdots=x_n$.

\section{GLOBAL CONSENSUS FOR MULTI-AGENT SYSTEMS WITH INPUT SATURATION}\label{secmain}
In this section, under the condition that the underlying graph is directed, we consider the multi-agent system subject to input saturation, i.e., system (\ref{system}), and with the distributed consensus control protocol (\ref{inputc}).  We show that the global consensus can be achieved . We first consider the case that the underlying graph is directed and strongly connected, then we consider the case the underlying graph is directed and has a spanning tree.
\subsection{Strongly Connected Digraphs}
In this subsection, we consider the situation that the underlying digraph is strongly connected, i.e., the Laplacian matrix $L$ is irreducible.
We have the following result.
\begin{theorem}\label{statictheorem}
Consider the multi-agent system (\ref{system}) with the distributed protocol (\ref{inputc}).
Suppose that the underlying graph $\mathcal G$ is directed and strongly connected.
Then global consensus is achieved.
\end{theorem}
\begin{proof}
Let $\xi$ be the vector defined in Lemma \ref{lemmaL}.
Let $x(t)=[x_1^\top(t),\dots,x_n^\top(t)]^\top$. Consider the following function:
\begin{align}\label{V}
V(x)=\sum_{i=1}^{n}\xi_i\sum_{l=1}^{p}\int_{0}^{-\sum_{j=1}^{n}L_{ij}c_l(x_{j}(t))}\sat_h(s)ds.
\end{align}
From Lemmas \ref{lemmaL} and \ref{lemma3}, we know $V(x)\ge0$ and from \ref{lemmaLx}, we know $V(x)=0$ if and only if $x_1=\cdots=x_n$.

The derivative of $V(x)$ along the trajectories of system (\ref{system}) with the distributed consensus protocol (\ref{inputc}) is
\begin{align}
&\dot{V}(x)=\frac{dV(x)}{dt}\nonumber\\
=&\sum_{i=1}^{n}\xi_i\sum_{l=1}^p[\sat_h(-\sum_{j=1}^{n}L_{ij}c_l(x_{j}(t)))][-\sum_{j=1}^{n}L_{ij}c_l(\dot{x}_{j}(t))]\nonumber\\
=&\sum_{i=1}^{n}\xi_i\sum_{l=1}^p[\sat_h(c_l(u_{i}(t)))][-\sum_{j=1}^{n}L_{ij}\sat_h(c_l(u_{j}(t)))]\nonumber\\
=&\sum_{i=1}^{n}\xi_i[\sat_h(u_i(t))]^\top\sum_{j=1}^{n}-L_{ij}\sat_h(u_j(t))\nonumber\\
=&-\sum_{i=1}^{n}\xi_iq_i(t),\label{dV}
\end{align}
where
\begin{align*}
q_{i}(t)=-\frac{1}{2}\sum_{j=1}^{n}L_{ij}\|\sat_h(u_{j}(t))-\sat_h(u_{i}(t))\|^2\ge0,
\end{align*}
and the last equality holds since
\begin{align}
&-\sum_{i=1}^{n}\xi_iq_i(t)\nonumber\\
=&\sum_{i=1}^{n}\frac{1}{2}\sum_{j=1}^{n}\xi_iL_{ij}\|\sat_h(u_{j}(t))-\sat_h(u_{i}(t))\|^2\nonumber\\
=&\sum_{i=1}^{n}\frac{1}{2}\sum_{j=1}^{n}\xi_iL_{ij}\Big[\|\sat_h(u_{j}(t))\|^2
+\|\sat_h(u_{i}(t))\|^2\Big]\nonumber\\
&-\sum_{i=1}^{n}\sum_{j=1}^{n}\xi_iL_{ij}[\sat_h(u_{j}(t))]^\top\sat_h(u_{i}(t))\nonumber\\
=&\frac{1}{2}\sum_{j=1}^{n}\|\sat_h(u_{j}(t))\|^2\sum_{i=1}^{n}\xi_iL_{ij}\nonumber\\
&+\frac{1}{2}\sum_{i=1}^{n}\xi_i\|\sat_h(u_{i}(t))\|^2\sum_{j=1}^{n}L_{ij}\nonumber\\
&-\sum_{i=1}^{n}\sum_{j=1}^{n}\xi_iL_{ij}[\sat_h(u_{j}(t))]^\top\sat_h(u_{i}(t))\nonumber\\
=&-\sum_{i=1}^{n}\sum_{j=1}^{n}\xi_iL_{ij}[\sat_h(u_{j}(t))]^\top\sat_h(u_{i}(t)),\label{xiqi}
\end{align}
where we have used $\xi^\top L=0$ and $L{\bf 1}_n=0$ in (\ref{xiqi}).

We then show that consensus is achieved and the input of each agent enters into the saturation level in finite time. From (\ref{dV}), we know that $\dot{V}=0$ if and only if $\sat_h(u_i(t))=\sat_h(u_j(t)),\forall i,j=1,\dots,n$. From Lemma \ref{lemmaLx}, this is equivalent to $x_i(t)=x_j(t),\forall i,j=1,\dots,n$. Thus  by LaSalle Invariance Principle \citep{khalil2002nonlinear}, we have
\begin{align}
\lim_{t\rightarrow\infty}x_i(t)-x_j(t)=0,i,j=1\dots,n,\label{stateir}
\end{align}
i.e., the consensus is achieved.

Finally, we estimate the convergence speed which will be used later.
Since $-\sum_{j=1}^{n}L_{ij}c_l(x_j(t)),i=1,\dots,n,l=1,\dots,p$ are continuous with respect to $t$, it then follows from (\ref{stateir})  that there exists a constant $T_1\ge0$ such that
\begin{align}
|c_l(u_i(t))|=\Big|-\sum_{j=1}^{n}L_{ij}c_l(x_j(t))\Big|\le h,\forall t\ge T_1.
\end{align}
In other words the saturation function in (\ref{system}) does not play a role after $T_1$.
Thus the multi-agent system (\ref{system}) with the distributed protocol (\ref{inputc}) reduces to
\begin{align}\label{systemws}
\dot{x}_i(t)=-\sum_{j=1}^{n}L_{ij}x_j(t),t\ge T_1.
\end{align}

Consider the following function
\begin{align}\label{Vtilde}
\tilde{V}(x)=\frac{1}{2}x^\top(t)(U\otimes I_p)x(t).
\end{align}
From Lemma \ref{lemmaU}, we know that $\tilde{V}(x(t))\ge0$.
The derivative of $\tilde{V}(x)$ along the trajectories of system (\ref{systemws}) satisfies
\begin{align*}
&\frac{d\tilde{V}(x)}{dt}=x^\top(t)(U\otimes I_p)\dot{x}(t)\nonumber\\
=&x^\top(t)(U\otimes I_p)(-L\otimes I_p)x(t)=-x^\top(t)(R\otimes I_p)x(t)\nonumber\\
\le&-\frac{\rho_2(R)}{\rho(U)}x^\top(t)(U\otimes I_p)x(t)=-2\frac{\rho_2(R)}{\rho(U)}\tilde{V}(x),\forall t\ge T_1.
\end{align*}
Thus
\begin{align}
\tilde{V}(x(t))\le \tilde{V}(x(T_1))e^{-2\frac{\rho_2(R)}{\rho(U)}(t-T_1)},\forall t\ge T_1.
\end{align}
Noting $\tilde{V}(x(t))$ is continuous with respect to $t$, there exists a positive constant $C_1$ such that
\begin{align*}
\tilde{V}(x(t))\le C_1,\forall t\in[0,T_1].
\end{align*}
Then
\begin{align}\label{tV}
\tilde{V}(x(t))\le C_2e^{-2\frac{\rho_2(R)}{\rho(U)}t},\forall t\ge 0,
\end{align}
where $C_2=\max\Big\{\tilde{V}(x(T_1)),C_1e^{2\frac{\rho_2(R)}{\rho(U)}T_1}\Big\}$.

Moreover, from Lemma \ref{lemmaU}, we know that
\begin{align}
&\sum_{j=1}^{n}\|u_j(t)\|^2=x^{\top}(t)(L^\top L\otimes I_p)x(t)\nonumber\\
\le&\frac{\rho(L^\top L)}{\rho_2(U)}x^{\top}(t)(U\otimes I_p)x(t)=2\frac{\rho(L^\top L)}{\rho_2(U)}\tilde{V}(x(t))\nonumber\\
\le&2\frac{\rho(L^\top L)}{\rho_2(U)}C_2e^{-2\frac{\rho_2(R)}{\rho(U)}t},\forall t\ge 0.\label{u}
\end{align}
Thus we complete the proof.
\end{proof}

\subsection{Digraphs Having A Spanning Tree}\label{secmainr}
In this subsection, we consider the case that the underlying graph is directed and has a spanning tree. In this case the corresponding Laplacian matrix $L$ is reducible. The following
mathematical methods are inspired by \citet{chen2007pinning} which are useful for our analysis. By proper
permutation, we can rewrite $L$ as the following Perron-Frobenius form:
\begin{eqnarray}
L=\left[\begin{array}{llll}L^{1,1}&L^{1,2}&\cdots&L^{1,M}\\
0&L^{2,2}&\cdots&L^{2,M}\\
\vdots&\vdots&\ddots&\vdots\\
0&0&\cdots&L^{M,M}
\end{array}\right],\label{PF}
\end{eqnarray}
where $L^{m,m}$ is with dimension $n_{m}$ and associated with the $m$-th
strongly connected component (SCC) of $\mathcal G$, denoted by $\scc_{m}$,
$m=1,\dots,M$.

Since $\mathcal G$ contains a spanning tree, then each $L^{m,m}$ is irreducible or
has one dimension and for each $m<M$, $L^{m,q}\ne 0$ for at least one
$q>m$. Define an auxiliary matrix
$\tilde{L}^{m,m}=[\tilde{L}^{m,m}_{ij}]_{i,j=1}^{n_{m}}$ as
\begin{eqnarray*}
\tilde{L}^{m,m}_{ij}=\begin{cases}L^{m,m}_{ij}&i\ne j,\\
-\sum_{r=1,r\not=i}^{n_{m}}L^{m,m}_{ir}&i=j.\end{cases}
\end{eqnarray*}

Let ${\xi^{m}}=[\xi^m_1,\dots,\xi^m_{n_m}]^\top$ be the positive left eigenvector of the irreducible
$\tilde{L}^{m,m}$ corresponding to the eigenvalue zero and has the sum of components equaling to $1$.
Denote $\Xi^{m}=diag[\xi^{m}]$, $Q^{m}=\frac{1}{2}[\Xi^{m}L^{m,m}+(\Xi^{m}L^{m,m})^{\top}]$, $m=1,\dots,M$, and $U^{M}=\Xi^{M}-\xi^{M}(\xi^{M})^{\top}$.
Then, we have
\begin{lemma}\label{lemmaQ}
Under the setup above, $Q^{m}$ is positive definite for all $m<M$. $Q^M$ and $U^M$ are semi-positive definite. Moreover,
\begin{align}\label{QU}
Q^M\ge \frac{\rho_2(Q^M)}{\rho(U^M)}U^M.
\end{align}
\end{lemma}
\begin{proof}
For the proof of $Q^{m}$ is positive definite for all $m<M$, please see \citet[Lemma 3.1]{wu2005synchronization}.
\end{proof}

Let $N_0=0,~N_{m}=\sum_{i=1}^{m}n_{i},m=1,\dots,M$. Then the $i$-th agent in $\scc_m$ is the $N_{m-1}+i$-th agent in the whole graph. In the following, we exchangeably use $v^m_i$ and $v_{N_{m-1}+i}$ to denote this agent. Accordingly, denote $x^m_i(t)=x_{N_{m-1}+i}(t)$, $\hat{x}^m_i(t)=\hat{x}_{N_{m-1}+i}(t)$, $u^m_i(t)=u_{N_{m-1}+i}(t)$ and define $u^{m}(t)=[(u_{1}^m)^{\top}(t),\dots,(u_{n_{m}}^m)^{\top}(t)]^{\top}$.

Our second main result result is given in the following theorem.

\begin{theorem}\label{dynamictheorem}
Consider the multi-agent system (\ref{system}) with the distributed consensus protocol (\ref{inputc}).
Suppose that the underlying graph $\mathcal G$ is directed and has a directed spanning tree, and $L$ is written in the form of (\ref{PF}).
Then global consensus is achieved.
\end{theorem}
We illustrate the main idea of the proof here. For the detail of the proof, please see Appendix \ref{appendixa}.
By Theorem \ref{statictheorem}, all agents in $\scc_M$ achieve consensus since $L^{M,M}$ is irreducible or zero. Then, all agents in $\scc_{M-1}$ which is either strongly connected or of dimension one achieve the same consensus value as those in $\scc_M$ since this is a leader following problem with agents in $\scc_M$ are leaders and agents in $\scc_{M-1}$ are followers. By applying the similar analysis, we see that  all agents in $\scc_{m}$, $m=M-2,\dots,1$, which is either strongly connected or of dimension one achieve the same consensus value as above since this is a leader following problem with agents in $\scc_M,\scc_{M_1},\dots,\scc_{m+1}$ are leaders and agents in $\scc_{m}$ are followers. Therefore, the result follows.

\begin{remark}
\citet{li2011consensus} also consider the consensus problem for multi-agent systems with input saturation constraint under the condition that the underlying graphes having directed spanning trees. The Lyapunov function used in this paper is different form the one used in \cite{li2011consensus}. It facilitates the design of event-triggered control laws as shown in Section \ref{secmaini}.
\end{remark}

\section{EVENT-TRIGGERED CONTROL FOR MULTI-AGENT SYSTEMS WITH INPUT SATURATION}\label{secmaini}

To implement consensus protocol (\ref{inputc}), continuous states from neighbours are needed. However, continuous communication is impractical in physical applications.
In order to avoid continuous sending of information among agents and updating of actuators, we equip the distributed consensus protocol (\ref{inputc}) with event-triggered communication scheme under which the control signal is only updated when the event-triggered condition is satisfied.
Here, we use the following multi-agent system with input saturation and event-triggered control protocol
\begin{align}
\dot{x}_i(t)&=\sat_h(\hat{u}_i(t)),~i=1,\dots,n, t\ge0,\label{systemi}\\
\hat{u}_i(t)&=-\sum_{j=1}^{n}L_{ij}x_j(t^j_{k_j(t)}),\label{inputi}
\end{align}
where $k_{j}(t)=\argmax_{k}\{t^{j}_{k}\le t\}$. The increasing time agent-wise sequence $\{t_{k}^{j}\}_{k=1}^{\infty}$, $j=1,\dots,n$, named {\em triggering time sequence} of agent $v_j$ which will be determined later. We also assume $t^j_1=0, j=1,\dots,n$.
Note that the control protocol (\ref{inputi}) only updates at the triggering times and is constant between consecutive triggering times.

For simplicity, let $\hat{x}_{i}(t)=x_{i}(t_{k_{i}(t)}^{i})$, 
$e_{i}(t)=\hat{x}_{i}(t)-x_{i}(t)$, and $f_{i}(t)=\sat_h(\hat{u}_{i}(t))-\sat_h(u_{i}(t))$. 

In the following, we show that global consensus is achieved for the multi-agent system (\ref{systemi}) with event-triggered control protocol (\ref{inputi}). Similar to the analysis in Section \ref{secmain}, we first consider the case that the underlying digraph is strongly connected, we then consider the case that the underlying digraph  has a spanning tree.

\subsection{Strongly Connected Digraphs}
In this subsection, we consider the situation that the
underlying graph is directed and strongly connected. We have the following result.
\begin{theorem}\label{statictheoreme}
Consider the multi-agent system (\ref{systemi}) with the even-triggered control protocol (\ref{inputi}).
Suppose that the underlying graph $\mathcal G$ is directed and strongly connected.
Given $\alpha_i>0,\beta_i>0$ and the first triggering time $t^i_1=0$, agent $v_i$ determines the triggering times $\{t^i_k\}_{k=2}^{\infty}$ by
\begin{align}\label{statictriggersingle}
t^i_{k+1}=\max_{r\ge t^i_k}\Big\{r:\|e_i(t)\|^2\le \alpha_ie^{-\beta_it},
\forall t\in[t^i_k,r]\Big\}.
\end{align}
Then (i) there is no Zeno behavior; (ii) global consensus is achieved.
\end{theorem}
\begin{proof}
(i) We prove that there is no Zeno behavior by contradiction. Suppose there exists Zeno behavior. Then there exists agent $v_i$, such that $\lim_{k\rightarrow\infty}t^i_k=T_0$ with $T_0$ is a constant. Let $\varepsilon_0=\frac{\sqrt{\alpha_i}}{2\sqrt{p}h}e^{-\frac{1}{2}\beta_iT_0}>0$. Then from the property of limit,  there exists a positive integer $N(\varepsilon_0)$ such that
\begin{align}\label{zeno}
t^i_k\in[T_0-\varepsilon_0,T_0],~\forall k\ge N(\varepsilon_0).
\end{align}
Noting $\|\sat_h(s)\|\le h\sqrt{p}$ for any $s\in\mathbb{R}^p$, we have $$\|\sat_h(\hat{u}_i(t))\|\le h\sqrt{p}.$$ Noting $$\bigg|\frac{d\|e_i(t)\|}{dt}\bigg|\le\|\dot{x}_i(t)\|=\|\sat_h(\hat{u}_i(t))\|\le h\sqrt{p},$$ and $\|\hat{x}_i(t^i_{k})-x_i(t^i_{k})\|=0$ for any triggering time $t^i_k$, we can conclude that one sufficient condition to guarantee the inequality in condition (\ref{statictriggersingle}) is
\begin{align}\label{suffi2}
(t-t^i_k)h\sqrt{p}\le\sqrt{\alpha_i}e^{-\frac{1}{2}\beta_it}.
\end{align}
Then
\begin{align}
&t^i_{N(\varepsilon_0)+1}-t^i_{N(\varepsilon_0)}
\ge\frac{\sqrt{\alpha_i}}{\sqrt{p}h}e^{-\frac{1}{2}\beta_it^i_{N(\varepsilon_0)+1}}\nonumber\\
&\ge\frac{\sqrt{\alpha_i}}{\sqrt{p}h}e^{-\frac{1}{2}\beta_iT_0}=2\varepsilon_0,
\end{align}
which contradicts to (\ref{zeno}). Therefore, there is no Zeno behavior.

(ii) Firstly, the derivative of $V(x)$ defined in (\ref{V}) along the trajectories of system (\ref{systemi}) with the even-triggered control protocol (\ref{inputi}) satisfies
\begin{align}
&\dot{V}=\frac{dV(x)}{dt}\nonumber\\
=&\sum_{i=1}^{n}\xi_i\sum_{l=1}^p[\sat_h(-\sum_{j=1}^{n}L_{ij}c_l(x_{j}(t)))][-\sum_{j=1}^{n}L_{ij}c_l(\dot{x}_{j}(t))]\nonumber\\
=&\sum_{i=1}^{n}\xi_i\sum_{l=1}^p[\sat_h(c_l(u_{i}(t)))][-\sum_{j=1}^{n}L_{ij}\sat_h(c_l(\hat{u}_{j}(t)))]\nonumber\\
=&-\sum_{i=1}^{n}\xi_i[\sat_h(u_i(t))]^\top\sum_{j=1}^{n}L_{ij}\sat_h(\hat{u}_j(t))\nonumber\\
=&-\sum_{i=1}^{n}\xi_i[\sat_h(\hat{u}_i(t))-f_i(t)]^\top\sum_{j=1}^{n}L_{ij}\sat_h(\hat{u}_j(t))\nonumber\\
=&-\sum_{i=1}^{n}\sum_{j=1}^{n}\xi_iL_{ij}[\sat_h(\hat{u}_i(t))]^\top\sat_h(\hat{u}_j(t))\nonumber\\
&+\sum_{i=1}^{n}\sum_{j=1}^{n}\xi_iL_{ij}[f_i(t)]^\top\sat_h(\hat{u}_j(t))\nonumber\\
\overset{*}{=}&\sum_{i=1}^{n}\frac{\xi_i}{2}\sum_{j=1}^{n}L_{ij}\|\sat_h(\hat{u}_i(t))-\sat_h(\hat{u}_j(t))\|^2\nonumber\\
&+\sum_{i=1}^{n}\sum_{j=1,j\neq i}^{n}\xi_iL_{ij}[f_i(t)]^{\top}[\sat_h(\hat{u}_j(t))-\sat_h(\hat{u}_i(t))]\nonumber\\
\le&\sum_{i=1}^{n}\frac{\xi_i}{2}\sum_{j=1}^{n}L_{ij}\|\sat_h(\hat{u}_j(t))-\sat_h(\hat{u}_i(t))\|^2\nonumber\\
&+\sum_{i=1}^{n}\sum_{j=1,j\neq i}^{n}\Big\{-\xi_iL_{ij}\frac{1}{4}\|\sat_h(\hat{u}_j(t))-\sat_h(\hat{u}_i(t))\|^2\nonumber\\
&-\xi_iL_{ij}\|f_i(t)\|^2\Big\}\nonumber\\
=&\sum_{i=1}^{n}\frac{\xi_i}{4}\sum_{j=1}^{n}L_{ij}\|\sat_h(\hat{u}_j(t))-\sat_h(\hat{u}_i(t))\|^2\nonumber\\
&+\sum_{i=1}^{n}\xi_iL_{ii}\|f_i(t)\|^2\nonumber\\
=&\sum_{i=1}^{n}\frac{\xi_i}{4}\sum_{j=1,j\neq i}^{n}L_{ij}\|\sat_h(u_j(t))+f_j(t)\nonumber\\
&-\sat_h(u_i(t))-f_i(t)\|^2+\sum_{i=1}^{n}\xi_iL_{ii}\|f_i(t)\|^2\nonumber\\
=&\sum_{i=1}^{n}\frac{\xi_i}{4}\sum_{j=1,j\neq i}^{n}L_{ij}\Big\{\|\sat_h(u_j(t))-\sat_h(u_i(t))\|^2\nonumber\\
&+2[\sat_h(u_j(t))-\sat_h(u_i(t))]^\top[f_j(t)-f_i(t)]\nonumber\\
&+\|f_j(t)-f_i(t)\|^2\Big\}+\sum_{i=1}^{n}\xi_iL_{ii}\|f_i(t)\|^2\nonumber\\
\le&\sum_{i=1}^{n}\frac{\xi_i}{4}\sum_{j=1,j\neq i}^{n}L_{ij}\Big\{\|\sat_h(u_j(t))-\sat_h(u_i(t))\|^2\nonumber\\
&-\frac{1}{2}\|\sat_h(u_j(t))-\sat_h(u_i(t))\|^2-2\|f_j(t)-f_i(t)\|^2\nonumber\\
&+\|f_j(t)-f_i(t)\|^2\Big\}+\sum_{i=1}^{n}\xi_iL_{ii}\|f_i(t)\|^2\nonumber\\
=&\sum_{i=1}^{n}\frac{\xi_i}{4}\sum_{j=1,j\neq i}^{n}L_{ij}\Big\{\frac{1}{2}\|\sat_h(u_j(t))-\sat_h(u_i(t))\|^2\nonumber\\
&-\|f_j(t)-f_i(t)\|^2\Big\}+\sum_{i=1}^{n}\xi_iL_{ii}\|f_i(t)\|^2\nonumber\\
\le&-\sum_{i=1}^{n}\frac{\xi_i}{4}q_i(t)+\sum_{i=1}^{n}\frac{\xi_i}{4}\sum_{j=1,j\neq i}^{n}L_{ij}\Big\{-2\|f_j(t)\|^2\nonumber\\
&-2\|f_i(t)\|^2\Big\}+\sum_{i=1}^{n}\xi_iL_{ii}\|f_i(t)\|^2\nonumber\\
=&-\sum_{i=1}^{n}\frac{\xi_i}{4}q_i(t)\nonumber\\
&+\sum_{i=1}^{n}\frac{\xi_i}{4}\Big\{\sum_{j=1}^{n}L_{ij}(-2)
\|f_j(t)\|^2+2L_{ii}\|f_i(t)\|^2\Big\}\nonumber\\
&-\sum_{i=1}^{n}\frac{\xi_i}{4}\sum_{j=1,j\neq i}^{n}L_{ij}2\|f_i(t)\|^2+\sum_{i=1}^{n}\xi_iL_{ii}\|f_i(t)\|^2\nonumber\\
=&-\sum_{i=1}^{n}\frac{\xi_i}{4}q_i(t)+\sum_{i=1}^{n}2\xi_iL_{ii}\|f_i(t)\|^2\nonumber\\
=&-\sum_{i=1}^{n}\frac{\xi_i}{4}q_i(t)+\sum_{i=1}^{n}2\xi_iL_{ii}\|\sat_h(\hat{u}_i(t))-\sat_h(u_i(t))\|^2\nonumber\\
\overset{**}{\le}&-\sum_{i=1}^{n}\frac{\xi_i}{4}q_i(t)+\sum_{i=1}^{n}2\xi_iL_{ii}\|\hat{u}_i(t)-u_i(t)\|^2\nonumber\\
=&-\sum_{i=1}^{n}\frac{\xi_i}{4}q_i(t)+\sum_{i=1}^{n}2\xi_iL_{ii}\Big\|\sum_{j=1}^nL_{ij}e_j(t)\Big\|^2\nonumber\\
\le&-\sum_{i=1}^{n}\frac{\xi_i}{4}q_i(t)+2\max_{i\in\mathcal I}\Big\{\xi_iL_{ii}\Big\}e^\top(t)(L^\top L\otimes I_p)e(t)\nonumber\\
\le&-\sum_{i=1}^{n}\frac{\xi_i}{4}q_i(t)+2\max_{i\in\mathcal I}\Big\{\xi_iL_{ii}\Big\}\rho(L^\top L)\sum_{i=1}^{n}\|e_i(t)\|^2,\label{dVi}
\end{align}
where the equality denoted by  $\overset{*}{=}$ holds due to (\ref{xiqi}) and the inequality denoted by $\overset{**}{\le}$ holds due to Lemma \ref{lemma3}.

We then show that global consensus is achieved and the input of each agent enters into the saturation level in finite time. Let's treat $z_i(t)=e^{-\beta_it},t\ge0$ as an additional state to agent $v_i,i\in\mathcal I$. And let $z=[z_1,\dots,z_n]^\top$ Consider a  Lyapunov candidate:
\begin{align*}
W(x,z)=V(x)+4\max\Big\{\xi_iL_{ii}\Big\}\rho(L^\top L)\sum_{i=1}^{n}\frac{\alpha_i}{\beta_i}z_i(t),
\end{align*}
where $V(x)$ is defined in (\ref{V}).
Then the derivative of $W(x,z)$ along the multi-agent system (\ref{systemi}) with the even-triggered control protocol (\ref{inputi}) and system $\dot{z}_i(t)=-\beta_iz_i(t)$ is
\begin{align}
&\dot{W}=\frac{dW(x,z)}{dt}\nonumber\\
=&\dot{V}(x)-4\max\Big\{\xi_iL_{ii}\Big\}\rho(L^\top L)\sum_{i=1}^{n}\alpha_ie^{-\beta_it}\nonumber\\
\le&-\sum_{i=1}^{n}\frac{\xi_i}{4}q_i(t)+2\max\Big\{\xi_iL_{ii}\Big\}\rho(L^\top L)\sum_{i=1}^{n}\|e_i(t)\|^2\nonumber\\
&-4\max\Big\{\xi_iL_{ii}\Big\}\rho(L^\top L)\sum_{i=1}^{n}\alpha_ie^{-\beta_it}\nonumber\\
\le&-\sum_{i=1}^{n}\frac{\xi_i}{4}q_i(t)-2\max\Big\{\xi_iL_{ii}\Big\}\rho(L^\top L)\sum_{i=1}^{n}\alpha_ie^{-\beta_it}\le0.\nonumber
\end{align}
Then by LaSalle Invariance Principle \citep{khalil2002nonlinear}, similar to the proof in Theorem \ref{statictheorem}, we have
\begin{align}\label{stateiri}
\lim_{t\rightarrow\infty}x_j(t)-x_i(t)=0,i,j\in\mathcal I.
\end{align}

Finally, we estimate the convergence speed which will be used later.
Since $c_l(\hat{u}_i(t))=-\sum_{j=1}^{n}L_{ij}c_l(x_j(t))-\sum_{j=1}^{n}L_{ij}c_l(e_j(t))$, (\ref{statictriggersingle}), $-\sum_{j=1}^{n}L_{ij}c_l(x_j(t)),i\in\mathcal I,l=1,\dots,p$ are continuous with respect to $t$, it then follows from (\ref{stateiri}) that there exists a constant $T_3\ge0$ such that
\begin{align}
|c_l(\hat{u}_i(t))|&\le\Big|-\sum_{j=1}^{n}L_{ij}c_l(x_j(t))\Big|+\Big|-\sum_{j=1}^{n}L_{ij}c_l(e_j(t))\Big|\nonumber\\
&\le h,\forall t\ge T_3.\label{uih}
\end{align}
In other words the saturation function in (\ref{systemi}) does not play a role after $T_3$.
Thus the multi-agent system (\ref{systemi}) with distributed protocol (\ref{inputi}) reduces to
\begin{align}\label{systemwsi}
\dot{x}_i(t)=-\sum_{j=1}^{n}L_{ij}\hat{x}_j(t),t\ge T_3.
\end{align}

Similar to the proof of \citet[Theorem 2]{yi2016distributed}, we can conclude that there exist $C_5>0$ and $C_6>0$ such that
\begin{align}
\tilde{V}(x(t))\le C_5e^{-C_6t},\forall t\ge T_3,
\end{align}
where $\tilde{V}(x)$ defined in (\ref{Vtilde}).
Similar to the way to get (\ref{tV}), we have
\begin{align}\label{tVi}
\tilde{V}(x(t))\le C_7e^{-C_6t},\forall t\ge 0,
\end{align}
where $C_7$ is a positive constant.

Moreover, similar to the analysis for obtaining (\ref{u}), we have
\begin{align}
\sum_{i=1}^{n}\|\hat{u}_i(t)\|^2=&\sum_{i=1}^{n}\|u_i(t)-\sum_{j=1}^nL_{ij}e_j(t)\|^2\nonumber\\
\le&2\sum_{i=1}^{n}\|u_i(t)\|^2+2\rho(L^\top L)\sum_{i=1}^{n}\|e_i(t)\|^2\nonumber\\
\le&C_7e^{-C_8t},\forall t\ge 0,\label{ui}
\end{align}
where $C_7$ and $C_8$ are two positive constants.
\end{proof}
\begin{remark}
We call the event-triggered control protocol (\ref{inputi}) together with the event-triggering condition (\ref{statictriggersingle}) a event-triggered control law. It is distributed since it only needs its own state information, without any priori knowledge of any global parameter, such as the eigenvalue of the Laplacian matrix.
\end{remark}


\subsection{Digraphs Having A Spanning Tree}
In this subsection, we consider the case that the underlying graph has a spanning tree. We use the same notations as in Section \ref{secmainr}. For simplicity, let $u^m_i(t)=u_{N_{m-1}+i}(t)$, $e^m_i(t)=e_{N_{m-1}+i}(t)$, $f^m_i(t)=f_{N_{m-1}+i}(t)$, $\alpha^m_i=\alpha_{N_{m-1}+i}$, $\beta^m_i=\beta_{N_{m-1}+i}$, and $u^{m}(t)=[(u_{1}^m)^{\top}(t),\dots,(u_{n_{m}}^m)^{\top}(t)]^{\top}$.

Our fourth main result is given in the following theorem.

\begin{theorem}\label{dynamictheoreme}
Consider the multi-agent system (\ref{systemi}) with the even-triggered control protocol (\ref{inputi}).
Suppose that the underlying graph $\mathcal G$ is directed and has a spanning tree, and $L$ is written in the form of (\ref{PF}). Given $\alpha_i>0,\beta_i>0$ and the first triggering time $t^i_1=0$, agent $v_i$ determines the triggering times $\{t^i_k\}_{k=2}^{\infty}$ by the event-triggering condition (\ref{statictriggersingle}).
Then (i) there is no Zeno behavior; (ii) global consensus is achieved.
\end{theorem}
The proof is similar to the proof of Theorem \ref{dynamictheorem}. For the details, please see Appendix \ref{appendixb}.

\section{SIMULATIONS}\label{secsimulation}

\begin{figure}[hbt]
\centering
\includegraphics[width=0.45\textwidth]{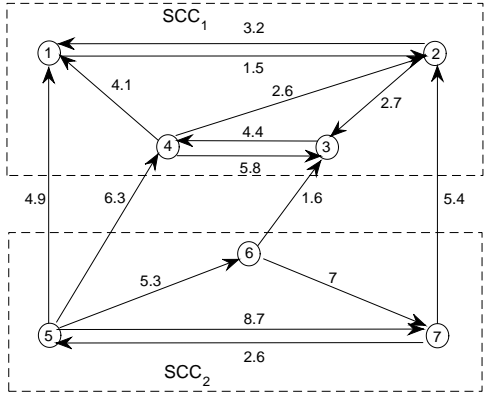}
\caption{The communication topology.}
\label{fig:1}
\end{figure}

\begin{figure}[hbt]
\begin{subfigure}{.5\textwidth}
  \centering
  \includegraphics[width=.9\linewidth]{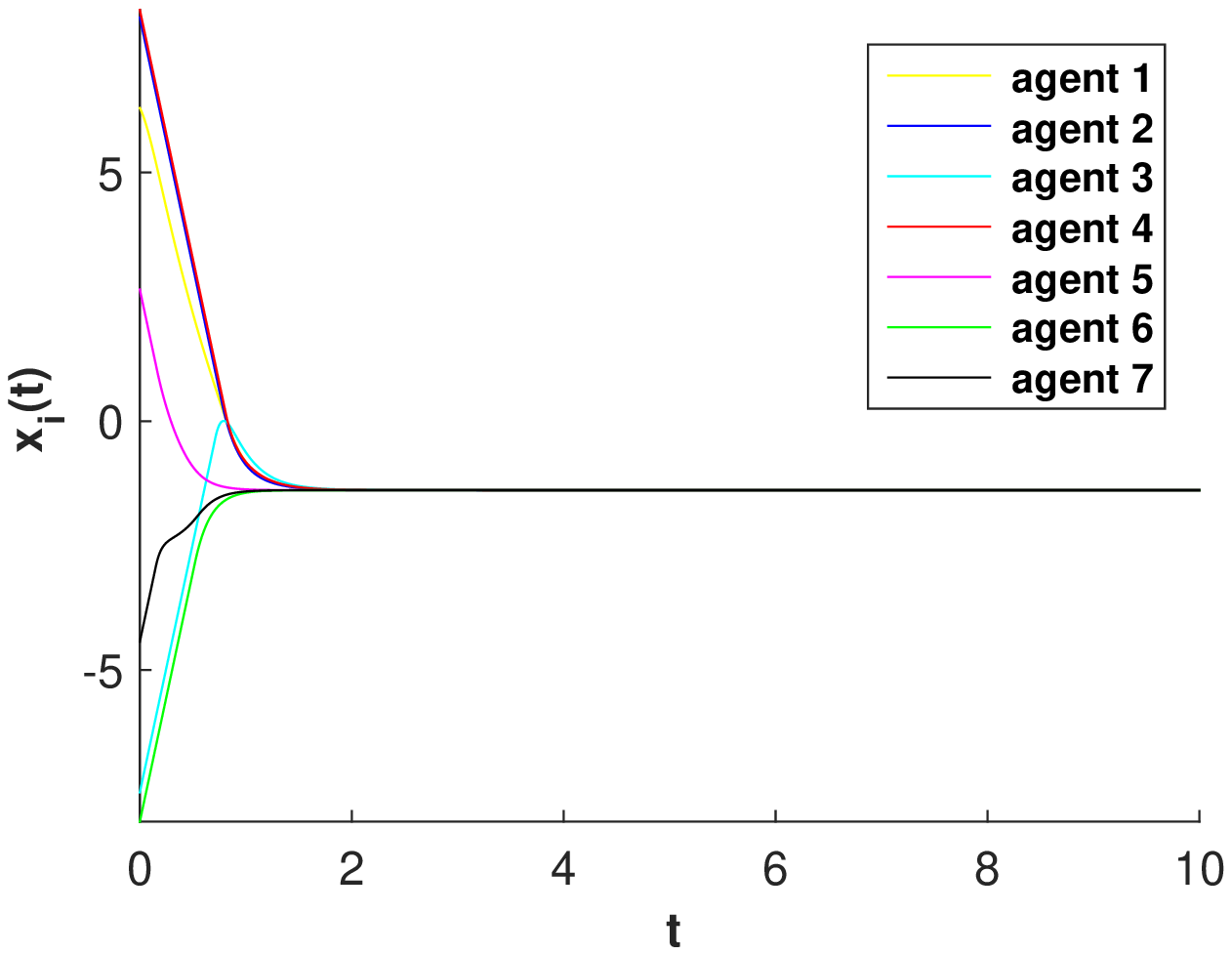}
  \caption{}
  \label{fig:2a}
\end{subfigure}%
\\
\begin{subfigure}{.5\textwidth}
  \centering
  \includegraphics[width=.9\linewidth]{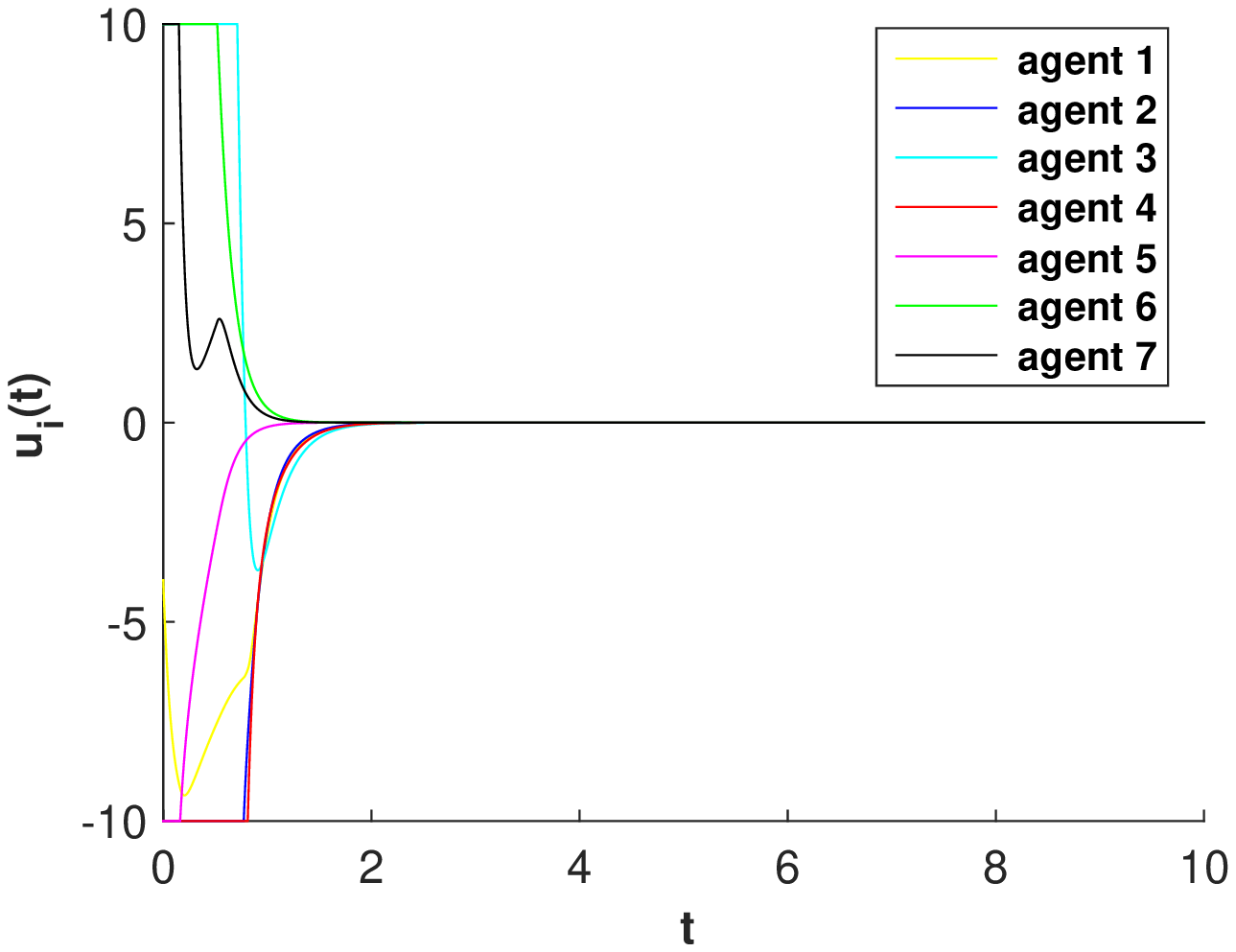}
  \caption{}
  \label{fig:2b}
\end{subfigure}
\caption{(a) The state evolution of the multi-agent system (\ref{system}) with the distributed protocol (\ref{inputc}). (b) The saturated input of each agent.}
\label{fig:2}
\end{figure}

\begin{figure}
\begin{subfigure}{.5\textwidth}
  \centering
  \includegraphics[width=.9\linewidth]{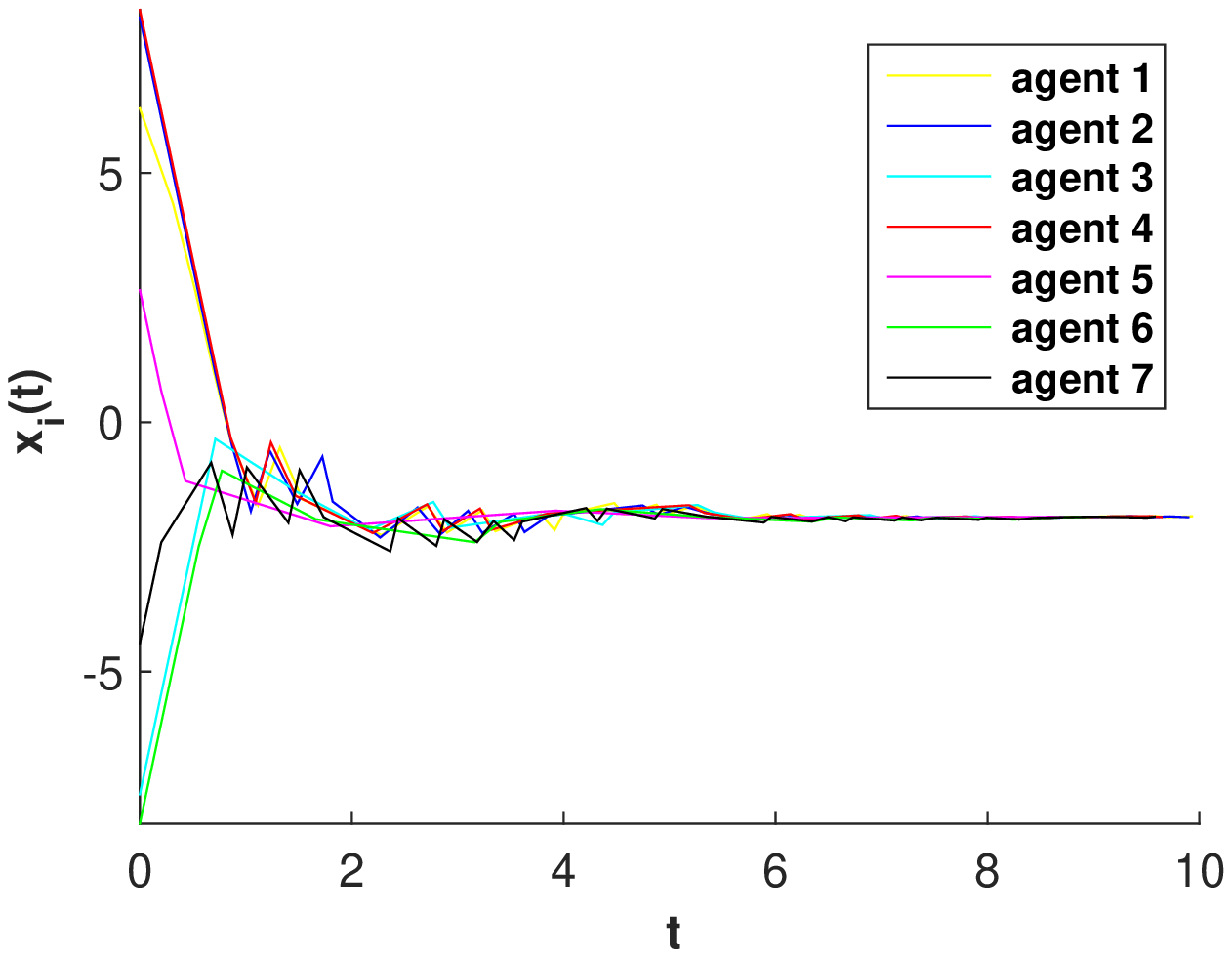}
  \caption{}
  \label{fig:3a}
\end{subfigure}%
\\
\begin{subfigure}{.5\textwidth}
  \centering
  \includegraphics[width=.9\linewidth]{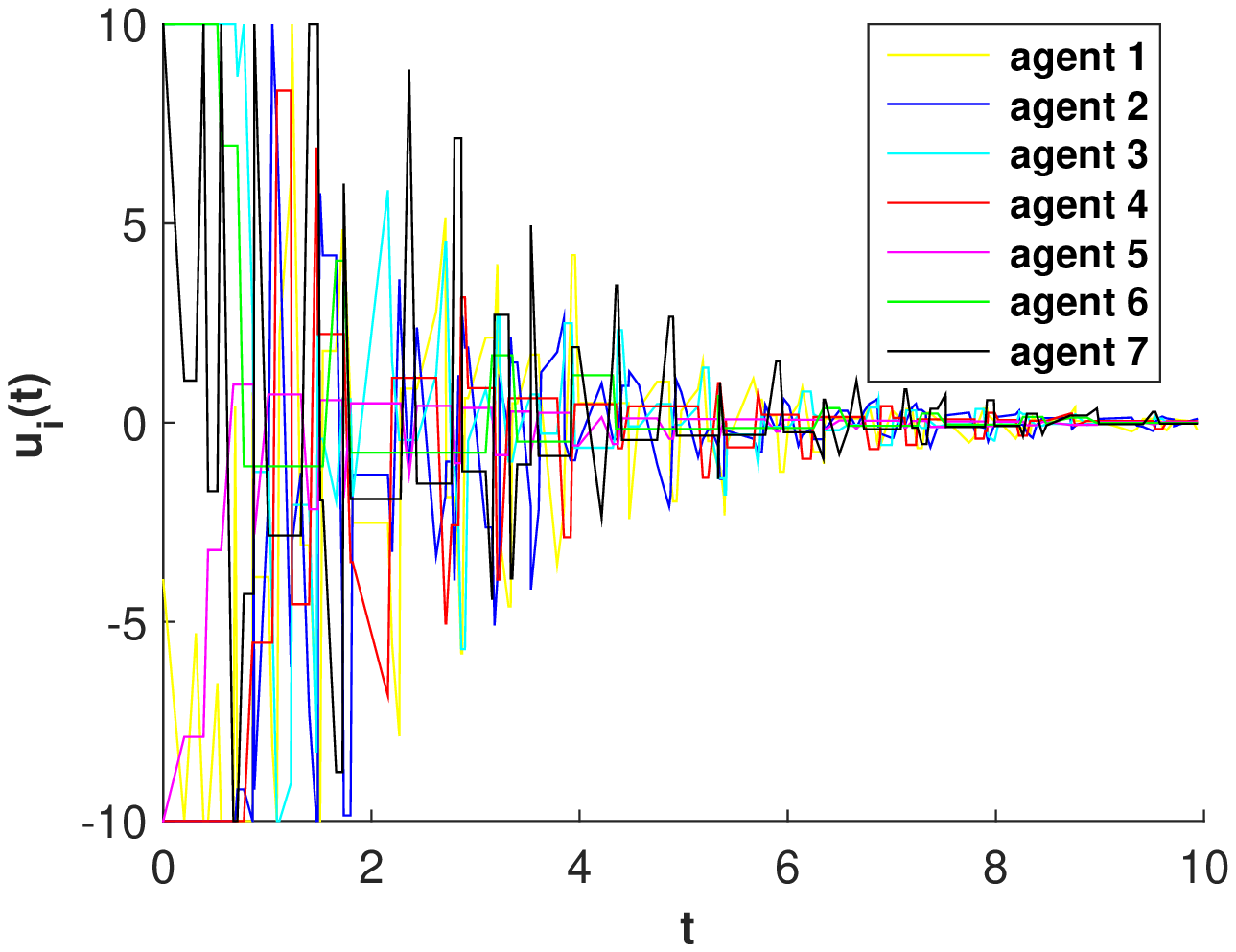}
  \caption{}
  \label{fig:3b}
\end{subfigure}
\\
\begin{subfigure}{.5\textwidth}
  \centering
  \includegraphics[width=.9\linewidth]{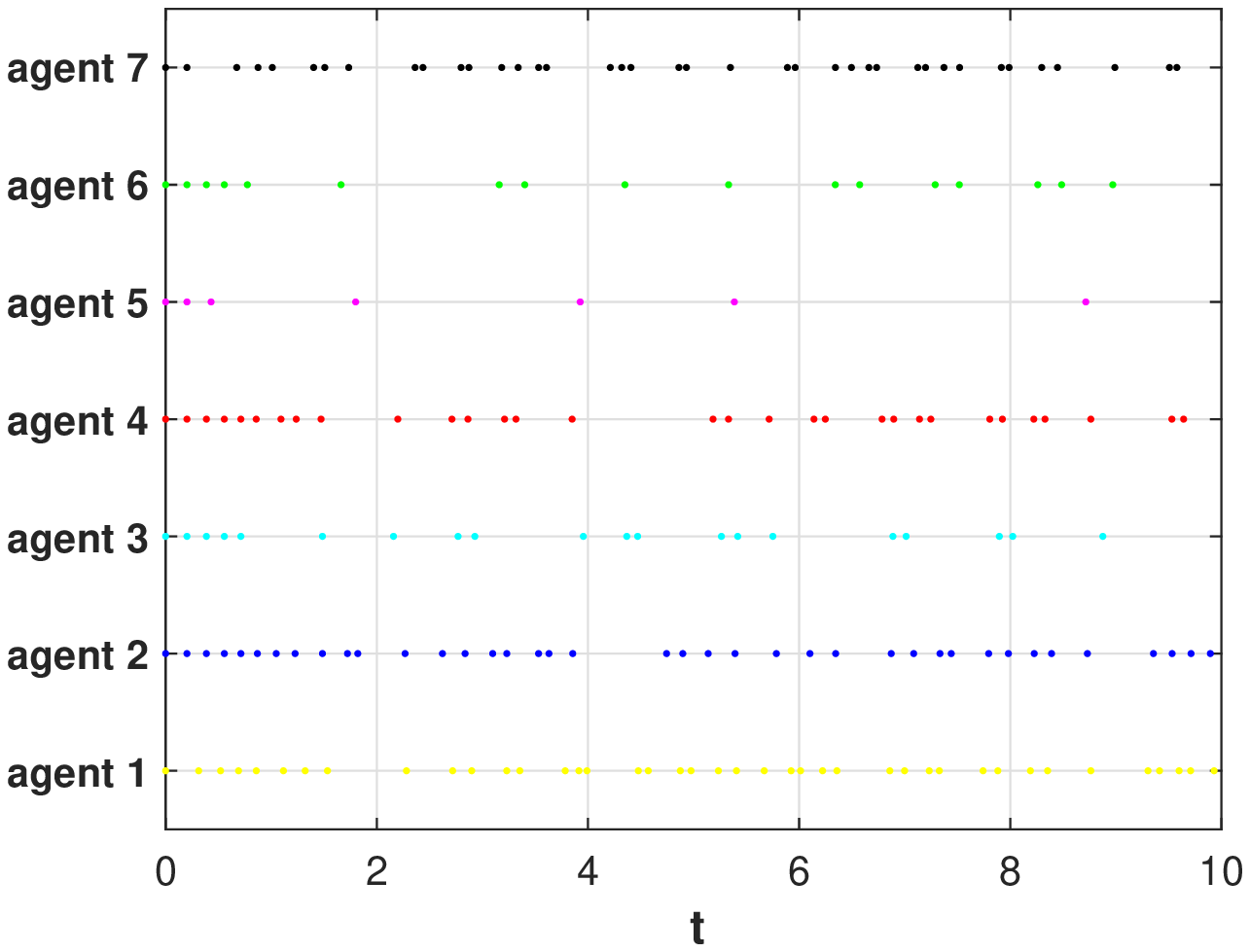}
  \caption{}
  \label{fig:3c}
\end{subfigure}
\caption{(a) The state evolution of the multi-agent system (\ref{systemi}) with the even-triggered control protocol (\ref{inputi}) under the event-triggering condition (\ref{statictriggersingle}). (b) The saturated input of each agent. (c) The triggering times for each agent.}
\label{fig:3}
\end{figure}

In this section, a numerical example is given to demonstrate the presented results. The saturation parameter is $h=10$.
Consider a directed graph of seven agents with the Laplacian matrix
\begin{eqnarray*}
L=\left[\begin{array}{rrrrrrr}12.2 &-3.2 & 0     &-4.1   &-4.9  &0    &0\\
    -1.5   &9.5 &0     &-2.6   &0    &0    &-5.4\\
    0     &-2.7  &10.1 &-5.8   &0    &-1.6  &0\\
    0     &0    &-4.4   &10.7 &-6.3  &0    &0\\
    0     &0    &0     &0     &2.6 &0    &-2.6\\
    0     &0    &0     &0     &-5.3  &5.3 &0\\
    0     &0    &0     &0     &-8.7  &-7   &15.7
\end{array}\right],
\end{eqnarray*}
whose topology is shown in Fig. \ref{fig:1}. The seven agents can be divided into two strongly connected components, i.e. the first four agents form a strongly connected component and the rest form another.
The initial value of each agent is randomly selected within the interval $[-10,10]$. Here, $x(0)=[6.2945,    8.1158,   -7.4603,    8.2675,    2.6472,   -8.0492, \\  -4.4300]^{\top}$. Fig. \ref{fig:2} (a) shows the state evolution of the multi-agent system (\ref{system}) with the distributed protocol (\ref{inputc}) and Fig. \ref{fig:2} (b) shows the saturated input of each agent. From Fig. \ref{fig:2} (a) and (b), we see that global consensus is achieved and $\sat_h(u_i(t))$ is within the saturation level.
Fig. \ref{fig:3} (a) shows the state evolution of the multi-agent system (\ref{systemi}) with the even-triggered control protocol (\ref{inputi}) under the event-triggering condition (\ref{statictriggersingle}) with $\alpha_i=10$ and $\beta_i=1$. Fig. \ref{fig:3} (b) shows the saturated input of each agent. Fig. \ref{fig:3} (c) shows the corresponding triggering times for each agent.
From Fig. \ref{fig:3} (a) and (b), we see that global consensus is achieved and $\sat_h(u_i(t))$ is within the saturation level. Moreover, from Fig. \ref{fig:3} (c), we see that each agent only needs to broadcast its state to its neighbors at its triggering times. Thus continuous broadcasting is avoided.

\section{CONCLUSION}\label{secconclusion}
In this paper, we studied the global consensus problem for multi-agent systems with input saturation constraints. The communication topologies among agents were described as directed graphs with directed spanning trees. We used a novel Lyapunov function to show that saturation constrain did not play a role after a finite time and global consensus is achieved. Moreover,
we presented a distributed event-triggered control law to reduce the  overall
need of communication and system updates. We showed global consensus is still achieved. Furthermore, the event-triggered control law was shown to be free of Zeno behavior.  Future research directions of this work include considering more general systems, time delays, and the self-triggered control.


\bibliography{refs}           

\begin{thebibliography}{33}
\expandafter\ifx\csname natexlab\endcsname\relax\def\natexlab#1{#1}\fi
\expandafter\ifx\csname url\endcsname\relax
  \def\url#1{\texttt{#1}}\fi
\expandafter\ifx\csname urlprefix\endcsname\relax\def\urlprefix{URL }\fi

\bibitem[{{\AA}str{\"o}m and Bernhardsson(1999)}]{aastrom1999comparison}
{\AA}str{\"o}m, K.~J., Bernhardsson, B., 1999. Comparison of periodic and event
  based sampling for first-order stochastic systems. In: Proceedings of the
  14th IFAC World congress. Vol.~11. Citeseer, pp. 301--306.

\bibitem[{Chen et~al.(2007)Chen, Liu, and Lu}]{chen2007pinning}
Chen, T., Liu, X., Lu, W., 2007. Pinning complex networks by a single
  controller. IEEE Transactions on Circuits and Systems I: Regular Papers
  54~(6), 1317--1326.

\bibitem[{Dimarogonas et~al.(2012)Dimarogonas, Frazzoli, and
  Johansson}]{dimarogonas2012distributed}
Dimarogonas, D.~V., Frazzoli, E., Johansson, K.~H., 2012. Distributed
  event-triggered control for multi-agent systems. IEEE Transactions on
  Automatic Control 57~(5), 1291--1297.

\bibitem[{Fan et~al.(2013)Fan, Feng, Wang, and Song}]{fan2013distributed}
Fan, Y., Feng, G., Wang, Y., Song, C., 2013. Distributed event-triggered
  control of multi-agent systems with combinational measurements. Automatica
  49~(2), 671--675.

\bibitem[{Heemels et~al.(2012)Heemels, Johansson, and
  Tabuada}]{heemels2012introduction}
Heemels, W., Johansson, K.~H., Tabuada, P., 2012. An introduction to
  event-triggered and self-triggered control. In: Decision and Control (CDC),
  2012 IEEE 51st Annual Conference on. IEEE, pp. 3270--3285.

\bibitem[{Horn and Johnson(2012)}]{horn2012matrix}
Horn, R.~A., Johnson, C.~R., 2012. Matrix analysis. Cambridge university press.

\bibitem[{Johansson et~al.(1999)Johansson, Egerstedt, Lygeros, and
  Sastry}]{johansson1999regularization}
Johansson, K.~H., Egerstedt, M., Lygeros, J., Sastry, S., 1999. On the
  regularization of {Z}eno hybrid automata. Systems \& Control Letters 38~(3),
  141--150.

\bibitem[{Khalil(2002)}]{khalil2002nonlinear}
Khalil, H.~K., 2002. Nonlinear systems, 3rd. Prentice-Hall, New Jersey.

\bibitem[{Kiener et~al.(2014)Kiener, Lehmann, and
  Johansson}]{kiener2014actuator}
Kiener, G.~A., Lehmann, D., Johansson, K.~H., 2014. Actuator saturation and
  anti-windup compensation in event-triggered control. Discrete event dynamic
  systems 24~(2), 173--197.

\bibitem[{Li et~al.(2011)Li, Xiang, and Wei}]{li2011consensus}
Li, Y., Xiang, J., Wei, W., 2011. Consensus problems for linear time-invariant
  multi-agent systems with saturation constraints. IET Control Theory \&
  Applications 5~(6), 823--829.

\bibitem[{Lim and Ahn(2016)}]{lim2016consensus}
Lim, Y.-H., Ahn, H.-S., 2016. Consensus with output saturations. arXiv preprint
  arXiv:1606.05980.

\bibitem[{Liu et~al.(2011)Liu, Lu, and Chen}]{Liu2011consensus}
Liu, B., Lu, W., Chen, T., 2011. Consensus in networks of multiagents with
  switching topologies modeled as adapted stochastic processes. SIAM Journal on
  Control and Optimization 49~(1), 227--253.

\bibitem[{Lu and Chen(2006)}]{lu2006new}
Lu, W., Chen, T., 2006. New approach to synchronization analysis of linearly
  coupled ordinary differential systems. Physica D: Nonlinear Phenomena
  213~(2), 214--230.

\bibitem[{Lu and Chen(2007)}]{lu2007new}
Lu, W., Chen, T., 2007. A new approach to synchronization analysis of linearly
  coupled map lattices. Chinese Annals of Mathematics, Series B 28~(2),
  149--160.

\bibitem[{Meng and Chen(2013)}]{meng2013event}
Meng, X., Chen, T., 2013. Event based agreement protocols for multi-agent
  networks. Automatica 49~(7), 2125--2132.

\bibitem[{Meng et~al.(2015)Meng, Xie, Soh, Nowzari, and
  Pappas}]{meng2015periodic}
Meng, X., Xie, L., Soh, Y.~C., Nowzari, C., Pappas, G.~J., 2015. Periodic
  event-triggered average consensus over directed graphs. In: Decision and
  Control (CDC), 2015 IEEE 54th Annual Conference on. IEEE, pp. 4151--4156.

\bibitem[{Meng et~al.(2013)Meng, Zhao, and Lin}]{meng2013global}
Meng, Z., Zhao, Z., Lin, Z., 2013. On global leader-following consensus of
  identical linear dynamic systems subject to actuator saturation. Systems \&
  Control Letters 62~(2), 132--142.

\bibitem[{Mesbahi and Egerstedt(2010)}]{mesbahi2010graph}
Mesbahi, M., Egerstedt, M., 2010. Graph Theoretic Methods in Multiagent
  Networks. Princeton University Press.

\bibitem[{Olfati-Saber and Murray(2004)}]{olfati2004consensus}
Olfati-Saber, R., Murray, R.~M., 2004. Consensus problems in networks of agents
  with switching topology and time-delays. IEEE Transactions on Automatic
  Control 49~(9), 1520--1533.

\bibitem[{Ren et~al.(2007)Ren, Beard, and Atkins}]{ren2007information}
Ren, W., Beard, R.~W., Atkins, E.~M., 2007. Information consensus in
  multivehicle cooperative control. IEEE Control Systems Magazine 27~(2),
  71--82.

\bibitem[{Seyboth et~al.(2013)Seyboth, Dimarogonas, and
  Johansson}]{seyboth2013event}
Seyboth, G.~S., Dimarogonas, D.~V., Johansson, K.~H., 2013. Event-based
  broadcasting for multi-agent average consensus. Automatica 49~(1), 245--252.

\bibitem[{Tabuada(2007)}]{tabuada2007event}
Tabuada, P., 2007. Event-triggered real-time scheduling of stabilizing control
  tasks. IEEE Transactions on Automatic Control 52~(9), 1680--1685.

\bibitem[{Wang and Sun(2016)}]{wang2016conditions}
Wang, Q., Sun, C., 2016. Conditions for consensus in directed networks of
  agents with heterogeneous output saturation. IET Control Theory \&
  Applications 10~(16), 2119--2127.

\bibitem[{Wang and Lemmon(2011)}]{wang2011event}
Wang, X., Lemmon, M.~D., 2011. Event-triggering in distributed networked
  control systems. IEEE Transactions on Automatic Control 56~(3), 586--601.

\bibitem[{Wu(2005)}]{wu2005synchronization}
Wu, C.~W., 2005. Synchronization in networks of nonlinear dynamical systems
  coupled via a directed graph. Nonlinearity 18~(3), 1057.

\bibitem[{Wu and Yang(2016)}]{wu2016distributed}
Wu, X., Yang, T., 2016. Distributed constrained event-triggered consensus: L 2
  gain design result. In: Industrial Electronics Society, IECON 2016-42nd
  Annual Conference of the IEEE. IEEE, pp. 5420--5425.

\bibitem[{Xie and Lin(2017)}]{xie2017event}
Xie, Y., Lin, Z., 2017. Event-triggered global stabilization of multiple
  integrator systems with bounded controls. In: Proceedings of the American
  Control Conference. IEEE.

\bibitem[{Yang et~al.(2016)Yang, Ren, Liu, and Chen}]{yang2016decentralized}
Yang, D., Ren, W., Liu, X., Chen, W., 2016. Decentralized event-triggered
  consensus for linear multi-agent systems under general directed graphs.
  Automatica 69, 242--249.

\bibitem[{Yang et~al.(2014)Yang, Meng, Dimarogonas, and
  Johansson}]{yang2014global}
Yang, T., Meng, Z., Dimarogonas, D.~V., Johansson, K.~H., 2014. Global
  consensus for discrete-time multi-agent systems with input saturation
  constraints. Automatica 50~(2), 499--506.

\bibitem[{Yi et~al.(2016{\natexlab{a}})Yi, Lu, and Chen}]{yi2016distributed}
Yi, X., Lu, W., Chen, T., 2016{\natexlab{a}}. Distributed event-triggered
  consensus for multi-agent systems with directed topologies. In: Proceedings
  of the 2016 Chinese Control and Decision Conference. IEEE, pp. 807--813.

\bibitem[{Yi et~al.(2017)Yi, Lu, and Chen}]{yi2017pull}
Yi, X., Lu, W., Chen, T., 2017. Pull-based distributed event-triggered
  consensus for multiagent systems with directed topologies. IEEE Transactions
  on Neural Networks and Learning Systems 28~(1), 71--79.

\bibitem[{Yi et~al.(2016{\natexlab{b}})Yi, Wei, Dimarogonas, and
  Johansson}]{yi2016formation}
Yi, X., Wei, J., Dimarogonas, D.~V., Johansson, K.~H., 2016{\natexlab{b}}.
  Formation control for multi-agent systems with connectivity preservation and
  event-triggered controllers. arXiv:1611.03105.

\bibitem[{You and Xie(2011)}]{you2011network}
You, K., Xie, L., 2011. Network topology and communication data rate for
  consensusability of discrete-time multi-agent systems. IEEE Transactions on
  Automatic Control 56~(10), 2262--2275.

\end{thebibliography}



\appendix
\section{Proof of Theorem \ref{dynamictheorem}}\label{appendixa}
For simplicity, hereby we only consider the case of $M=2$. The case $M>2$ can be treated in the similar manner.

Firstly, let's consider the second strongly connected component. All agents in $\scc_2$ do not dependent on any agents in $\scc_1$. Thus, the second strongly connected component can be treated as a strongly connected directed graph. Then form Theorem \ref{statictheorem}, we have
\begin{align*}
\lim_{t\rightarrow+\infty}x_i^2(t)-x_j^2(t)=0,i,j=1\dots,n_2.
\end{align*}
Then, from $-\sum_{j=1}^{n_2}c_l(x^2_{i}(t)),i=1,\dots,n_2,l=1,\dots,p$ are continuous with respect to $t$, we can conclude that there exists a constant $T_2\ge0$ such that
\begin{align}\label{state}
|c_l(u^2_{i}(t))|=\Big|-\sum_{j=1}^{n_2}L_{ij}^{2,2}c_l(x^2_{j}(t))\Big|\le h,\forall t\ge T_2.
\end{align}
In addition, similar to (\ref{u}), we have
\begin{align}
\|u^2(t)\|^2=\sum_{j=1}^{n_2}\|u^2_j(t)\|^2\le C_3e^{-C_4t}, t\ge0,
\end{align}
where $C_3$ and $C_4$ are two positive constants.

Secondly, let's consider the first strongly connected component.
Similar to $V(x)$ defined in (\ref{V}), define
\begin{align}
V_1(x)&=\sum_{i=1}^{n_1}\xi^1_i\sum_{l=1}^{p}\int_{0}^{c_l(u^1_i(t))}\sat_h(s)ds,\label{V1}\\
V_2(x)&=\sum_{i=1}^{n_2}\xi^2_i\sum_{l=1}^{p}\int_{0}^{c_l(u^2_{i}(t))}\sat_h(s)ds.\label{V2}
\end{align}
From the definition of the component operator $c_l(\cdot)$, we know $c_l(u^1_i(t))=-\sum_{j=1}^{n_1}L_{ij}^{1,1}c_l(x^1_{i}(t))
-\sum_{j=1}^{n_2}L_{ij}^{1,2}c_l(x^2_{i}(t))$ and $c_l(u^2_i(t))=
-\sum_{j=1}^{n_2}L_{ij}^{2,2}c_l(x^2_{i}(t))$. From Lemma \ref{lemma3},
we know $V_1(x)\ge0$ and $V_2(x)\ge0$.

Similar to (\ref{dV}), the derivative of $V_2(x)$ along the trajectories of system (\ref{system}) with the distributed consensus control protocol (\ref{inputc}) is
\begin{align}
\dot{V}_2(x)=\frac{dV_2(x)}{dt}=\sum_{i=1}^{n_2}-\xi_i^2q^2_i(t),
\end{align}
where
\begin{align*}
q^2_{i}(t)=-\frac{1}{2}\sum_{j=1}^{n}L_{ij}^{2,2}\|\sat_h(u^2_{j}(t))-\sat_h(u^2_{i}(t))\|^2\ge0.
\end{align*}
Moreover, similar to the proof of Theorem \ref{statictheorem}, we know $\dot{V}_2=0$ if and only if $x^2_i(t)=x^2_j(t),\forall i,j=1,\dots,n_2$

The derivative of $V_1(x)$  along the trajectories of system (\ref{system}) with the distributed consensus control protocol (\ref{inputc}) satisfies
\begin{align}
&\frac{dV_1(x)}{dt}=\sum_{i=1}^{n_1}\xi^1_i\sum_{l=1}^p\sat_h(c_l(u^1_i(t)))c_l(\dot{u}^1_i(t))\nonumber\\
=&\sum_{i=1}^{n_1}\xi^1_i\sum_{l=1}^pc_l(\sat_h(u^1_i(t)))\Big[-\sum_{j=1}^{n_1}L_{ij}^{1,1}c_l(\sat_h(u^1_j(t)))\nonumber\\
&-\sum_{j=1}^{n_2}L_{ij}^{1,2}c_l(\sat_h(u^2_j(t)))\Big]\nonumber\\
=&\sum_{i=1}^{n_1}\xi^1_i[\sat_h(u^1_i(t))]^\top\Big[-\sum_{j=1}^{n_1}L_{ij}^{1,1}\sat_h(u^1_j(t))\nonumber\\
&-\sum_{j=1}^{n_2}L_{ij}^{1,2}\sat_h(u^2_j(t))\Big]\nonumber\\
=&-[\sat_h(u^1(t))]^\top(Q^1\otimes I_p)\sat_h(u^1(t))\nonumber\\
&-\sum_{i=1}^{n_1}\xi^1_i[\sat_h(u^1_i(t))]^\top\sum_{j=1}^{n_2}L_{ij}^{1,2}\sat_h(u^2_j(t))\nonumber\\
\le&-\rho_2(Q^1)\|\sat_h(u^1(t))\|^2+\frac{\rho_2(Q^1)}{2}\sum_{i=1}^{n_1}\|\sat_h(u^1_i(t))\|^2\nonumber\\
&+\frac{1}{2\rho_2(Q^1)}\sum_{i=1}^{n_1}\Big\|\xi^1_i\sum_{j=1}^{n_2}L_{ij}^{1,2}\sat_h(u^2_j(t))\Big\|^2\nonumber\\
\le&-\frac{\rho_2(Q^1)}{2}\|\sat_h(u^1(t))\|^2\nonumber\\
&+\frac{n_1n_2\max\{(L_{ij}^{1,2})^2\}}{2\rho_2(Q^1)}\|\sat_h(u^2(t))\|^2\nonumber\\
\le&-\frac{\rho_2(Q^1)}{2}\|\sat_h(u^1(t))\|^2\nonumber\\
&+\frac{n_1n_2\max\{(L_{ij}^{1,2})^2\}}{2\rho_2(Q^1)}C_3e^{-C_4t}, t\ge 0.
\end{align}

Let's treat $y_i(t)=e^{-C_4t},t\ge0,i\in\mathcal I$ as an additional state of each agent. And let $y=[y_1,\dots,y_n]^\top$

Consider a Lyapunov candidate:
\begin{align}
V_3(x,y)=&V_1(x)+V_2(x)\nonumber\\
&+\frac{2n_1n_2\max\{(L_{ij}^{1,2})^2\}}{2\rho_2(Q^1)C_4n}C_3\sum_{i=1}^ny_i(t).
\end{align}

The derivative of $V_3(t)$ along the trajectories of system (\ref{system}) with distributed consensus control protocol (\ref{inputc}) is
\begin{align*}
\frac{dV_3(x,y)}{dt}=&\dot{V}_1(x)+\dot{V}_2(x)\\
&-\frac{2n_1n_2\max\{(L_{ij}^{1,2})^2\}}{2\rho_2(Q^1)n}C_3\sum_{i=1}^ny_i(t).
\end{align*}
Then,  we have
\begin{align}
\frac{dV_3(x,y)}{dt}
\le&-\frac{\rho_2(Q^1)}{2}\|\sat_h(u^1(t))\|^2+\sum_{i=1}^{n_2}-\xi_i^2q^2_i(t)\nonumber\\
&-\frac{n_1n_2\max\{(L_{ij}^{1,2})^2\}}{2\rho_2(Q^1)n}C_3\sum_{i=1}^ny_i(t),t\ge0.
\end{align}
Then by LaSalle Invariance Principle \citep{khalil2002nonlinear}, similar to the proof in Theorem \ref{statictheorem}, we have
\begin{align*}
\lim_{t\rightarrow\infty}x_j(t)-x_i(t)=0,\forall i,j\in\mathcal I.
\end{align*}
Similar to the proof in Theorem \ref{statictheorem}, we can show that after a finite time the saturation would not play a role any more, and global consensus is achieved.

\section{Proof of Theorem \ref{dynamictheoreme}}\label{appendixb}
(i) The proof of excluding Zeno behavior is the same as the proof of its counterpart in Theorem \ref{statictheoreme}. Thus we omit it here.

(ii) For simplicity, hereby we only consider the case of $M=2$. The case $M>2$ can be treated in the similar manner.

Firstly, let's consider the second strongly connected component. All agents in $\scc_2$ do not dependent on any agents in $\scc_1$. Thus, the second strongly connected component can be treated as a strongly connected directed graph. Then form Theorem \ref{statictheoreme}, we have
\begin{align*}
\lim_{t\rightarrow\infty}x_i^2(t)-x_j^2(t)=0,i,j=1\dots,n_2.
\end{align*}
and there exists a constant $T_4\ge0$ such that
\begin{align}\label{stateire}
|c_l(\hat{u}^2_{i}(t))|=\Big|-\sum_{j=1}^{n_2}L_{ij}^{2,2}c_l(\hat{x}^2_{j}(t))\Big|\le h,\forall t\ge T_4.
\end{align}
In addition, similar to (\ref{ui}), we have
\begin{align}
\|\hat{u}^2(t)\|^2=\sum_{j=1}^{n_2}\|\hat{u}^2_j(t)\|^2\le C_9e^{-C_{10}t}, t\ge0,
\end{align}
where $C_9$ and $C_{10}$ are two positive constants.

Secondly, let's consider the first strongly connected component.
Similar to (\ref{dVi}), the derivative of $V_2(x)$ defined in (\ref{V2}) along  the trajectories of system (\ref{systemi}) with the event-triggered control protocol (\ref{inputi}) satisfies
\begin{align}
&\dot{V}_2(x)=\frac{dV_2(x)}{dt}\nonumber\\
\le&-\sum_{i=1}^{n_2}\frac{\xi_i^2}{4}q^2_i(t)+d_1\sum_{i=1}^{n_2}\|e^2_i(t)\|^2,
\end{align}
where
$$d_1=2\max_{i\in\mathcal I}\Big\{\xi^2_iL^{2,2}_{ii}\Big\}\rho((L^{2,2})^\top L^{2,2}).$$

The derivative of $V_1(x)$ defined in (\ref{V1}) along  the trajectories of system (\ref{systemi}) with the event-triggered control protocol (\ref{inputi}) satisfies
\begin{align}
&\dot{V}_1(x)=\frac{dV_1(x)}{dt}=\sum_{i=1}^{n_1}\xi^1_i\sum_{l=1}^p\sat_h(c_l(u^1_i(t)))c_l(\dot{u}^1_i(t))\nonumber\\
=&\sum_{i=1}^{n_1}\xi^1_i\sum_{l=1}^pc_l(\sat_h(u^1_i(t)))\Big[-\sum_{j=1}^{n_1}L_{ij}^{1,1}c_l(\sat_h(\hat{u}^1_j(t)))\nonumber\\
&-\sum_{j=1}^{n_2}L_{ij}^{1,2}c_l(\sat_h(\hat{u}^2_j(t)))\Big]\nonumber\\
=&\sum_{i=1}^{n_1}\xi^1_i[\sat_h(u^1_i(t))]^\top\Big[-\sum_{j=1}^{n_1}L_{ij}^{1,1}\sat_h(\hat{u}^1_j(t))\nonumber\\
&-\sum_{j=1}^{n_2}L_{ij}^{1,2}\sat_h(\hat{u}^2_j(t))\Big]\nonumber\\
=&\sum_{i=1}^{n_1}\xi^1_i[\sat_h(\hat{u}^1_i(t))-f^1_i(t)]^\top\Big[-\sum_{j=1}^{n_1}L_{ij}^{1,1}\sat_h(\hat{u}^1_j(t))\nonumber\\
&-\sum_{j=1}^{n_2}L_{ij}^{1,2}\sat_h(\hat{u}^2_j(t))\Big]\nonumber\\
=&-[\sat_h(\hat{u}^1(t))]^\top(Q^1\otimes I_p)\sat_h(\hat{u}^1(t))\nonumber\\
&+\sum_{i=1}^{n_1}\xi^1_i[\sat_h(\hat{u}^1_i(t))]^\top\sum_{j=1}^{n_2}L_{ij}^{1,2}\sat_h(\hat{u}^2_j(t))\nonumber\\
&+\sum_{i=1}^{n_1}\xi^1_i[f^1_i(t)]^\top\Big[\sum_{j=1}^{n_1}L_{ij}^{1,1}\sat_h(\hat{u}^1_j(t))\nonumber\\
&+\sum_{j=1}^{n_2}L_{ij}^{1,2}\sat_h(\hat{u}^2_j(t))\Big]\nonumber\\
\le&-\rho_2(Q^1)\|\sat_h(\hat{u}^1(t))\|^2+\frac{\rho_2(Q^1)}{4}\sum_{i=1}^{n_1}\|\sat_h(\hat{u}^1_i(t))\|^2\nonumber\\
&+\frac{1}{\rho_2(Q^1)}\sum_{i=1}^{n_1}\Big\|\xi^1_i\sum_{j=1}^{n_2}L_{ij}^{1,2}\sat_h(\hat{u}^2_j(t))\Big\|^2\nonumber\\
&+\frac{\rho_2(Q^1)}{4}\sum_{j=1}^{n_1}\|\sat_h(\hat{u}^1_j(t))\|^2\nonumber\\
&+\frac{1}{\rho_2(Q^1)}\sum_{j=1}^{n_1}\Big\|\sum_{i=1}^{n_1}\xi^1_iL_{ij}^{1,1}f^1_i(t)\Big\|^2\nonumber\\
&+\sum_{i=1}^{n_1}\frac{1}{4}\|f^1_i(t)\|^2+\sum_{i=1}^{n_1}\Big\|\xi^1_i\sum_{j=1}^{n_2}L_{ij}^{1,2}\sat_h(\hat{u}^2_j(t))\Big\|^2\nonumber\\
\le&-\frac{\rho_2(Q^1)}{2}\|\sat_h(\hat{u}^1(t))\|^2+d_2\sum_{i=1}^{n_1}\|f^1_i(t)\|^2\nonumber\\
&+d_3\|\sat_h(\hat{u}^2(t))\|^2,\label{dVme}
\end{align}
where
\begin{align*}
d_2=&\frac{1}{4}+(n_1)^2\max_{i\in\{1,\dots,n_1\}}\{(\xi^1_iL_{ij}^{1,1})^2\}\frac{1}{\rho_2(Q^1)},\\
d_3=&2n_1n_2\max_{i\in\{1,\dots,n_1\}}\{(\xi^1_iL_{ij}^{1,2})^2\}\Big(\frac{1}{\rho_2(Q^1)}+1\Big).
\end{align*}
Similar to the analysis to get (\ref{dVi}), from (\ref{dVme}), we have
\begin{align}
\dot{V}_1\le&-\frac{\rho_2(Q^1)}{2}\|\sat_h(\hat{u}^1(t))\|^2+d_4\sum_{i=1}^{n_1}\|e^1_i(t)\|^2\nonumber\\
&+d_4\sum_{i=1}^{n_2}\|e^2_i(t)\|^2+d_3\|\sat_h(\hat{u}^2(t))\|^2,
\end{align}
where
\begin{align*}
d_4=d_2\rho(L^\top L).
\end{align*}

Let's treat $\eta^r_i(t)=e^{-\beta^r_{i}y},t\ge0$ as an additional state of agent $v^r_i,r=1,2,i=1,\dots,n_2$, and $\theta^2_i(t)=e^{-C_{10}t},t\ge0$ as an additional state of agent $v^2_i,i=1,\dots,n_2$, and $\theta^1_i(t)=0,t\ge0$ as an additional state of agent $v^1_i,i=1,\dots,n_1$. For simplicity. let $\eta=[\eta^1_1,\dots,\eta^1_{n_1},\eta^2_1,\dots,\eta^1_{n_2}]^\top$ and $\theta=[\theta^1_1,\dots,\theta^1_{n_1},\theta^2_1,\dots,\theta^1_{n_2}]^\top$.

Consider the following Lyapunov candidate:
\begin{align}
W_r(x,\eta,\theta)=&V_1(x)+V_2(x)+2\frac{C_9}{C_{10}}d_3\sum_{i=1}^{n_2}\theta^2_i\nonumber\\
&+2\sum_{i=1}^{n_2}\frac{(d_1+d_4)\alpha^2_i}{\beta^2_i}\eta^2_i+2\sum_{i=1}^{n_1}\frac{d_4\alpha^1_i}{\beta^1_i}\eta^1_i.
\end{align}
The derivative of $W_r(t)$ along  the trajectories of system (\ref{systemi}) with the event-triggered control protocol (\ref{inputi}) satisfies
\begin{align*}
\frac{dW_r(x,\eta,\theta)}{dt}=&\dot{V}_1(x)+\dot{V}_2(x)-2C_9d_3\sum_{i=1}^{n_2}\theta^2_i\\
&-2\sum_{i=1}^{n_2}(d_1+d_4)\alpha^2_i\eta^2_i-2\sum_{i=1}^{n_1}d_4\alpha^1_i\eta^1_i.
\end{align*}
Then, for any $t\ge T_4$, we have
\begin{align}
&\frac{dW_r(x,\eta,\theta)}{dt}
\le-\frac{\rho_2(Q^1)}{2}\|\sat_h(u^1(t))\|^2+\sum_{i=1}^{n_2}-\frac{\xi_i^2}{4}q^2_i(t)\nonumber\\
&-C_9d_3\sum_{i=1}^{n_2}\theta^2_i-\sum_{i=1}^{n_2}(d_1+d_4)\alpha^2_i\eta^2_i-\sum_{i=1}^{n_1}d_4\alpha^1_i\eta^1_i.
\end{align}
Then by LaSalle Invariance Principle \citep{khalil2002nonlinear}, similar to the proof in Theorem \ref{statictheorem}, we have
\begin{align*}
\lim_{t\rightarrow\infty}x_j(t)-x_i(t)=0,i,j\in\mathcal I.
\end{align*}
Similar to the proof in Theorem \ref{statictheorem}, we can show that after a finite time the saturation would not play a role any more, and global consensus is achieved.

\end{document}